\documentclass[11pt]{article}
\usepackage{amsmath,amsfonts,latexsym,amssymb,chicago,annals3,graphicx,pictexwd}

\newcommand{\mycite}[1]{{\small \sc \citeNP{#1}}}

\def\R{{\mathbb R}}

\def\E{{\mathbb E}}

\def\labda1{\lambda_1}
\def\labda2{\lambda_2}

\def\comment#1{\relax}

\def\=in{\mathop{\rm =}}

\newtheorem{lemma}{Lemma}[section]

\newcommand{\eps}{{\varepsilon }}

\renewcommand{\mycite}[1]{{\small \sc \citeNP{#1}}}

\newtheorem{thm}{Theorem}[section]
\newtheorem{cor}[thm]{Corollary}
\newtheorem{defi}[thm]{Definition}

\newcommand{\la}{{\lambda}}

\renewcommand{\R}{{\mathbb R}}

\begin{document}

\title{Behavior of a second class particle in Hammersley's process}
\author{Eric Cator and Sergei Dobrynin}
\date{\today}
\affiliation{Delft University of Technology} \AMSsubject{Primary: 60C05,60K35, secondary 60F05.}
 \maketitle

\begin{abstract}
In the case of a rarefaction fan in a non-stationary Hammersley process, we explicitly calculate
the asymptotic behavior of the process as we move out along a ray, and the asymptotic distribution
of the angle within the rarefaction fan of a second class particle and a dual second class
particle. Furthermore, we consider a stationary Hammersley process and use the previous results to
show that trajectories of a second class particle and a dual second class particles touch with
probability one, and we give some information on the area enclosed by the two trajectories, up
until the first intersection point. This is linked to the area of influence of an added Poisson
point in the plane.
\end{abstract}

\section{Introduction}

\begin{figure}[!h]
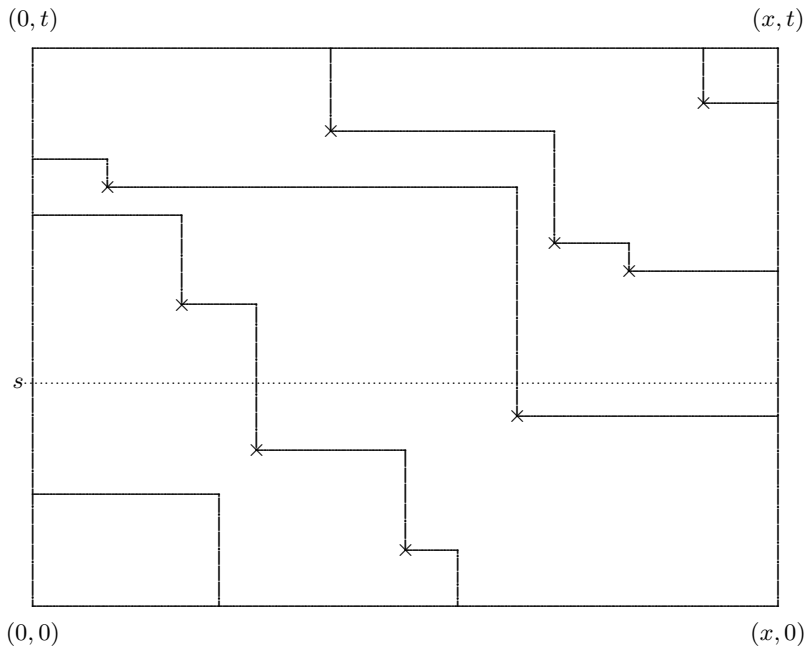

\begin{center}
\strut
\beginpicture
\footnotesize
\setcoordinatesystem units <0.06\textwidth,0.045\textwidth>
  \setplotarea x from 0 to 10, y from 0 to 12
\multiput {\bf $\times$} at
5 1
3 2.8
2 5.4
1 7.5
6.5 3.4
4 8.5
8 6
7 6.5
9 9
/
%

\put {$(0,0)$} at 0.0 -0.5
\put {$(x,0)$} at 10 -0.5
\put {$(0,t)$} at 0.0 10.5
\put {$(x,t)$} at 10 10.5
\put {$s$} at -0.2 4

\setlinear 
\plot 2.5 0.0  2.5 2  0.0 2 
/
\plot 5.7 0  5.7 1  5 1  5 2.8  3 2.8   
3 5.4  2 5.4  2 7 0 7
/
\plot 10 3.4  6.5 3.4  6.5 7.5  1 7.5  1 8  0 8
/
\plot 10 6   8 6  8 6.5  7 6.5  7 8.5  4 8.5 4 10
/
\plot 10 9  9 9  9 10
/
\plot 0 10  10 10  10 0 0 0 0 10
/
\setdots <2pt>
\plot -0.1 4 10 4
/
\endpicture

\caption{Space-time paths of the Hammersley's process, with sources
and sinks.} \label{fig:Hammsourcesink}
\end{center}
\end{figure}

In \mycite{ham:72}, a discrete interacting particle process is introduced to study the behavior of
the length of longest increasing subsequences of random permutations. In \mycite{aldi:95}, this
discrete process is generalized to a continuous time interacting particle process on the real line,
and they use the ergodic decomposition theorem to show local convergence to a Poisson process, when
moving out along a ray. In this paper, we will consider Hammersley's process with sources and
sinks, as introduced in \mycite{gro:02}. For an extensive description of this process, we refer to
\mycite{CaGr}, since our results will be partly based on results derived in that paper. Here we
will suffice with a brief description, based on Figure \ref{fig:Hammsourcesink}.

We consider the space-time paths of particles that started on the $x$-axis as sources, distributed
according to a Poisson distribution and we consider the $t$-axis as a time axis. In the positive
quadrant we have a Poisson process of what we call $\alpha$-points (denoted in Figure
\ref{fig:Hammsourcesink} by $\times$). At the time an $\alpha$-point appears, the particle
immediately to the right of it jumps to the location of the $\alpha$-point. Finally, we have a
Poisson process of sinks on the $t$-axis. Each sink makes the leftmost particle disappear. All
three Poisson processes are assumed to be independent. To know the particle configuration at time
$s$, we intersect a line at time $s$ with the space-time paths.

In \mycite{CaGr}, a connection was made between the continuous time Hammersley process and the
behavior of second class particles, which are well studied in the literature on discrete
interacting particle systems such as TASEP; see for example \mycite{lig:99}. For the Hammersley
process, it is natural to consider two types of second class particles: the usual one, where one
adds an extra particle at the origin, and a dual second class particle, which corresponds to adding
an extra sink at the origin (or removing the leftmost particle). In fact, the trajectories of these
two particles correspond to the two longest paths of the time-reversed process such that all
possible longest paths fall between these two longest paths. We can study the trajectories of these
two particles at the same time, that is, for one realization of the Hammersley process. In this
paper we study the behavior of a second class particle and its dual particle in the case of a
rarefaction fan, a phenomenon often observed in interacting particle systems. In \mycite{FerKip}
this problem is considered for TASEP. In \mycite{timo:02} the Hammersley process is considered with
general initial conditions, but the rarefaction fan is not treated. Also, our methods are quite
different and build more on the ideas of \mycite{CaGr}. In the final section we will study the
interaction between a second class particle and its dual in the case of a stationary Hammersley
process, and show that they will touch with probability one. We also study the area between the two
trajectories up until this point of touch. We did not find results in the literature on discrete
interacting particle systems that were similar to the results of our last section, so this
interaction phenomenon may be a specific feature of the Hammersley process.

\section{Second class particles in a rarefaction fan}

Let $\la, \mu$ be two positive reals, such that $\la\mu<1$. Let $t\mapsto L_{\la,\mu}(\cdot,
t)$ be Hammersley's process developing in time $t$, generated by a Poisson process of sources on
the positive x-axis of intensity $\la$, a Poisson process of sinks on the time axis of intensity
$\mu$ and a Poisson process on $\R^2_+$ of intensity $1$, where these Poisson processes are
independent. Here, $L_{\la,\mu}(\cdot, t)$ signifies the counting process that counts the number of
Hammersley particles on the half-line $(0,\infty)\times \{t\}$. As was shown in \mycite{gro:02},
the case $\la \mu =1$ corresponds to a stationary Hammersley process, which means that for each
$t\geq 0$, $L_{\la,\mu}(\cdot ,t)$ is a Poisson process with intensity $\la$.

As we mentioned in the introduction, we will consider two kinds of second class particles. A
``normal'' second class particle is created by putting an extra source in the origin. A {\em dual
second class particle} is created by putting an extra sink in the origin. The trajectories
$(X_t,t)$ of a second class particle and $(X'_t,t)$ of a dual second class particle are shown in
Figure \ref{fig:X&X'}. Note that we always have $X'_t\geq X_t$.

\begin{figure}[!ht]
\begin{center}
\strut
\beginpicture
\footnotesize
\setcoordinatesystem units < 8mm, 6mm>
  \setplotarea x from 0 to 12, y from 0 to 10.5
  \axis bottom  
    shiftedto y=0
    /
  \axis left 
    shiftedto x=0.0
    /
\multiput {\bf $\times$} at
5 1
3 2.8
2 5.4
1 7.5
6.5 3.4
4 8.5
8 6
7 6.5
11.5 5
11 2
/
    

\put {$t$} at -0.5 10
\put {$x$} at 9.1 -0.5
\put {$\left(X(t),t\right)$} at 7 10.5
\put {$\left(X(t)',t\right)$} at 11.5 10.5
\put {$(x,t)$} at 9.1 10.5

\setlinear 
\plot 2.5 0.0  2.5 2  0.0 2 
/
\plot 5.7 0  5.7 1  5 1  5 2.8  3 2.8   
3 5.4  2 5.4  2 7 0 7
/
\plot 10 0 10 3.4  6.5 3.4  6.5 7.5  1 7.5  1 8  0 8
/
\plot 12 2 11 2 11 6   8 6  8 6.5  7 6.5  7 8.5  4 8.5 4 10
/
\plot 11.5 10 11.5 5 12 5
/

{\setplotsymbol ({$*$}) 
\plotsymbolspacing=10pt
\plot 0 0 2.5 0 2.5 2 5 2 5 2.8 10 2.8 10 3.4 11 3.4 11 6 11.5 6 11.5 10
/

\setplotsymbol ({$\diamond$}) 
\plot 0 0 0 2  2.5 2  2.5 5.4  3 5.4  3 7.5  6.5 7.5  6.5 8 6.5 8.5 7 8.5  7 10
/}

\setdots <2pt>

\plot 0 10 12 10
/

\plot 9 0 9 10
/

\endpicture

\caption{Trajectories of $\left(X_t,t\right)$ and
$\left(X'_t,t\right)$} \label{fig:X&X'}
\end{center}
\end{figure}

Now consider the reversed process, where we use the North exits through $[0,x]\times \{t\}$ as
sources, the East exits through $\{x\}\times [0,t]$ as sinks and the $\beta$-points (these are the
upper-right corners of the space-time paths) as our Poisson process in $[0,x]\times [0,t]$. Burke's
Theorem for Hammersley's process (see \mycite{CaGr}) shows that this process is again a stationary
Hammersley process, if we start with a stationary process. It is not hard to see from Figure
\ref{fig:X&X'} that the trajectories of $X$ and $X'$ correspond to longest paths from $(x,t)$ to
$(0,0)$ in the reversed process.

A second class particle and a dual second class particle are symmetrical with respect to the main
diagonal in the following sense: the reflected trajectory of a dual second class particle is a
trajectory of a "normal" second class particle for the reflected realization, that is, with the
sources and sinks intensity exchanged.

A very important property of second class and dual second class particles is given in the following lemma, which is slightly more general than Lemma 2.2 in \mycite{CaGr}, but the proofs of these two lemma's are very similar, and therefore omitted here.

\begin{lemma}\label{lem:leaveout}
Consider a Hammersley process $L$ with some configuration of sources and sinks. Define two coupled Hammersley processes $\overline{L}$ and $\underline{L}$ such that all three processes use the same $\alpha$-points. Furthermore, $\overline{L}$ uses the same sources as $L$, but only a subset of the sinks of $L$, whereas $\underline{L}$ uses the same sinks as $L$, but only a subset of the sources of $L$. Denote with $X$ and $X'$ a second class particle respectively a dual second class particle for the process $L$. Then the space-time paths of $L$ and $\overline{L}$ coincide below the trajectory of $X$, whereas the space-time paths of $L$ and $\underline{L}$ coincide above the trajectory of $X'$.
\end{lemma}

The reader can verify this lemma by looking at Figure \ref{fig:X&X'}.

The results in this section are also valid when $\la$ and/or $\mu$ are equal to 0, if we keep in
mind that whenever $\la=0$, a dual second class particle remains on the $x$-axis (so $X'_t=+\infty$
for $t>0$), whereas when $\mu=0$, a second class particle remains on the $t$-axis (so $X_t=0$ for
$t>0$).

For $\rho >0$, we denote by $t\mapsto L_\rho (\cdot, t)$ a stationary Hammersley process, generated
by a process of sources of intensity $\rho $, a process of sinks of intensity $1/\rho $ and a
process of $\alpha$-points of intensity 1.

Throughout this section we are interested in the asymptotic behavior of the non-stationary process
$L_{\la,\mu}$ and trajectories of its second class and dual second class particles. We recall here
the corresponding results for a stationary process, obtained in \mycite{CaGr}.

\begin{thm}[Cator and Groeneboom]\label{thm:catgr} Let $t\rightarrow L_{\rho }(\cdot,t)$
be a stationary Hammersley process with intensity of sources $\rho $ and intensity of sinks $1/\rho
$. Let $X_t$ be the $x$-coordinate of a second class particle at time $t$ and $X'_t$ be the
$x$-coordinate of a dual second class particle at time $t$.

\noindent Then:

\noindent (i) $t^{-1} X_t \rightarrow 1/{\rho }^2$ almost surely as $t\rightarrow \infty$

\noindent (ii) $t^{-1} X'_t \rightarrow 1/{\rho }^2$ almost surely as $t\rightarrow \infty$
\end{thm}

This theorem combines the results of Theorem 2.1 and Remark 2.1 of \mycite{CaGr}. We can now show
the following theorem, describing the asymptotic local intensities of the process $L_{\la,\mu}$
when moving out along a ray $t=ax$, for $a>0$.

\begin{thm}\label{thm:1} Let $a>0$. Consider the random
particle configuration with counting process
$$
y \mapsto L_{\la,\mu}(x+y,ax)- L_{\la,\mu}(x,ax), y\geq -x.
$$

\noindent Then:

 \noindent(i) If $a>1/\mu^2$, the process converges in
 distribution, as $x\rightarrow \infty$, to a homogeneous Poisson
 process on $\R$, with intensity $1/\mu$.

 \noindent(ii) If $1/\mu^2>a>\la^2$, the process converges in
 distribution, as $x\rightarrow \infty$, to a homogeneous Poisson
 process on $\R$, with intensity $\sqrt{a}$.

\noindent(iii) If $\la^2>a$, the process converges in
 distribution, as $x\rightarrow \infty$, to a homogeneous Poisson
 process on $\R$, with intensity $\la$.
\end{thm}
\begin{figure}[!ht]
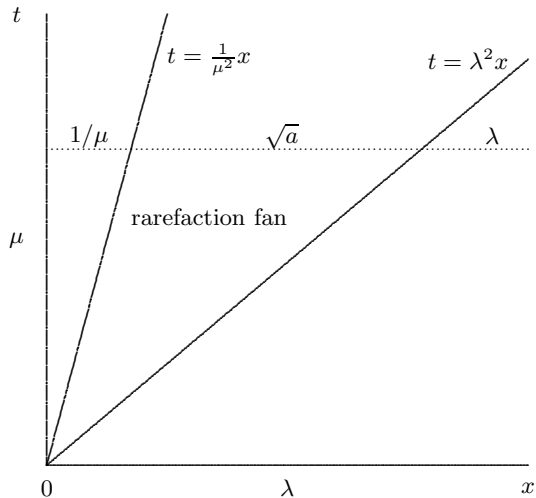

\begin{center}
\strut
\beginpicture
\footnotesize
\setcoordinatesystem units < 8mm, 6mm>
  \setplotarea x from 0 to 10, y from 0 to 12

\put {$0$} at 0 -0.5
\put {$\la$} at 4 -0.5
\put {$\mu$} at -0.5 5
\put {$t = {1\over{\mu}^2}x$} at 2.7 9
\put {$t = {\la}^2x$} at 7 9
\put {$x$} at 8 -0.5
\put {$1/\mu$} at 0.7 7.3
\put {$\la$} at 7.4 7.3
\put {$\sqrt{a}$} at 3.9 7.3
\put {$t$} at -0.5 10
\put {rarefaction fan} at 2.7 5.5

\setlinear 

\plot 0 10 0 0 8 0
/
\plot 0 0 8 9
/
\plot 0 0 2 10
/

\setdots <2pt>
\plot 0.0 7 8 7
/

\endpicture

\caption{Illustration of Theorem \ref{thm:1}} \label{fig:thm1}
\end{center}
\end{figure}

\noindent {\bf Proof:}\\
{\bf (i)} We couple the process $L_{\la,\mu}$ with a stationary process $L_{1/\mu}$ via the same
realization of the Poisson process in $\R^2_+$, the same realization of sinks on the $t$-axis, and
a "thickened" set of sources. This means that we independently add to the Poisson process of
intensity $\la$ on the $x$-axis a Poisson process of intensity $1/\mu-\lambda$. This way of
coupling (when only one set is changed (thinned or thickened) and the two others are kept fixed)
will be used throughout the paper without further explanation.

Let $(X'_t,t)$ be the trajectory of an isolated dual second class particle at zero of the process
$L_{1/\mu}$. As the process $L_{\la,\mu}$ can be obtained from $L_{1/\mu}$ by removing the added
sources, the space-time paths of both processes coincide above the trajectory $(X'_t,t)$ (Lemma \ref{lem:leaveout}). Now it is
enough to apply Theorem \ref{thm:catgr}, stating that $X'_t/t$ almost surely converges to
$\mu^2$, and to use the stationarity of $L_{1/\mu}$.\\
{\bf (ii)} This proof is analogous to the proof of Corollary 2.1, (i), in \mycite{CaGr}. Let $a'$
be such that $1/\mu^2>a'>a$. Note that $\sqrt{a'}>\la$ and $1/\sqrt{a'}>\mu$. Now we couple the
process $L_{\la,\mu}$ with a process $L_{\sqrt{a'}}$ generated by the same realization of points on
the plane and thickened sets of sources and sinks of intensities $\sqrt{a'}$ and $1/\sqrt{a'}$
correspondingly.

Below the trajectory $(X_t,t)$ of an isolated second class particle at zero of $L_{\sqrt{a'}}$, the
space-time paths of $L_{\sqrt{a'}}$ coincide with the space-time paths of the process $\bar
L_{\sqrt{a'}}$, which is obtained by removing all sinks (again Lemma \ref{lem:leaveout}).

Note that for any $t$, the set of intersections of space-time paths of the process $L_{\la,\mu}$
with the horizontal level $t$ is contained in the set of intersection of space-times paths of $\bar
L_{\sqrt{a'}}$ with the horizontal level $t$. In other words, at all times the process $\bar
L_{\sqrt{a'}}$ contains the particles of $L_{\la,\mu}$, and may have more particles, since the
first process has more sources and less sinks. It follows that the set of intersections of
space-time paths of $L_{\la,\mu}$ with any horizontal interval lying to the right of the trajectory
$(X_t,t)$ is contained in the set of intersections of paths of $L_{\sqrt{a'}}$ with this interval.

Similarly, we consider $a''$ such that $a>a''>\la^2$ and couple $L_{\la,\mu}$ with
$L_{\sqrt{a''}}$.  Now we have that {\em above} the trajectory $(X'_t,t)$ of an isolated dual
second class particle at zero of $L_{\sqrt{a''}}$, the space-time paths of $L_{\sqrt{a''}}$
coincide with the paths of ${\underline L}_{\sqrt{a''}}$, which is obtained by removing all
sources. This time we conclude that the set of intersections of space-time paths of $L_{\la,\mu}$
with any horizontal interval lying to the left of the trajectory $(X'_t,t)$ {\em contains} the set
of intersections of paths of $L_{\sqrt{a''}}$ with this interval.

The situation is now as follows: if we fix $R>0$ and consider the
process $L_{\la,\mu}$ on the horizontal interval $[x-R,x+R]\times
\{ax\}$, we know that for  any fixed $a'$ and $a''$ such that
$1/\mu^2>a'>a>a''>\la^2$ and for $x$ big enough, the particles of
this process are a subset of the stationary process
$L_{\sqrt{a'}}$ and a superset of the stationary process
$L_{\sqrt{a''}}$.
 Now we are exactly in the situation of the proof of the Corollary 2.1, (i), \mycite{CaGr} and
 statement (ii) follows by considering the Laplace transform of the number of particles in fixed
 disjoint intervals and using dominated convergence.\\
{\bf (iii)} The proof of the last statement proceeds along the same lines as the proof of statement
(i). Couple $L_{\la,\mu}$ with $L_{\la}$, note that their space-time paths coincide below the
trajectory $(X_t,t)$ of an isolated second class particle at zero of the process $L_{\la}$ and use
Theorem \ref{thm:catgr} and the stationarity of $L_{\la}$.\\
\mbox{ } \hfill $\Box$\\

Let $F_{\la,\mu}(x,t)$ be the flux of the process $L_{\la,\mu}$ at $(x,t)$, which we define by the
number of space-time paths of $L_{\la,\mu}$ contained in the rectangle $[0,x]\times [0,t]$. Note
that $F_{\la,\mu}(x,t)$ can be represented as the number of crossings of the interval
$[0,x]\times\{t\}$ by space-time paths of $L_{\la,\mu}$ plus the number of crossings of the
interval $\{0\}\times [0,t]$. It is also equal to the number of crossings of the interval
$[0,x]\times\{0\}$ plus the number of crossings of the interval $\{x\}\times [0,t]$ by space-time
paths of $L_{\la,\mu}$.

\begin{defi}
When moving out along a ray $t=ax$, we define the local intensity function $I(x,h)$ by the
following formula:
\[I(x,h)=\E(F_{\la,\mu}(x+h,ax)-F_{\la,\mu}(x,ax)).\]
We suppress the dependence on $a, \la$ and $\mu$ in our notation.
\end{defi}

Note that
$$
I(x,h)= \E(L_{\la,\mu}(x+h,ax)-L_{\la,\mu}(x,ax)),
$$
so it is equal to the expected number of crossings of the interval $[x,x+h]\times\{ax\}$ by
space-time paths of the process $L_{\la,\mu}$.

\begin{cor}\label{cor:pwl}
Fix the positive constants $a, \la$ and $\mu$. We have
$$
\lim\limits_{x\rightarrow \infty} I(x,h) = \begin{cases} h/\mu, & a\geq 1/{\mu^2}, \\
\sqrt{a}\cdot h, & 1/{\mu^2}>a\geq \la^2 \\ \la\cdot h,  & \la^2>a\geq 0.\end{cases}
$$
\end{cor}

\noindent {\bf Proof:} Coupling in a standard way the process $L_{\la,\mu}$ with the process
$L_{1/\mu}$ (adding sources), we can estimate
$$
L_{\la,\mu}(x+h,ax)-L_{\la,\mu}(x,ax)\leq L_{1/\mu}(x+h,ax)-L_{1/\mu}(x,ax)
$$
for any realization and any $x\geq 0$. Define for $x\geq 0$
\[ A_x=L_{\la,\mu}(x+h,ax)-L_{\la,\mu}(x,ax) \mbox{  and  } B_x=L_{1/\mu}(x+h,ax)-L_{1/\mu}(x,ax).\]
Define a random variable $A$, whose distribution depends on the fixed parameters $a, \la$ and
$\mu$, as follows:
\[ A\sim \begin{cases} {\rm Pois}(h/\mu), & a\geq 1/{\mu^2}, \\
{\rm Pois}(\sqrt{a}\cdot h), & 1/{\mu^2}>a\geq \la^2 \\ {\rm Pois}(\la\cdot h),  & \la^2>a\geq
0.\end{cases} \] Theorem \ref{thm:1} shows that
\[ A_x\stackrel{\rm d}{\longrightarrow} A\ \ (x\to \infty).\]
Furthermore,
\[ \forall\ x\geq 0:\ B_x\sim {\rm Pois}(h/\mu).\]
For any $M>0$ we have
\[ \E (A_x) = \E (A_x\wedge M) + \E ((A_x-M)1_{\{A_x>M\}}).\]
Also,
\[ \forall\ x>0:\ \E ((A_x-M)1_{\{A_x>M\}})\leq \E ((B_x-M)1_{\{B_x>M\}})=\E ((B_0-M)1_{\{B_0>M\}}).\]
Since
\[ \lim_{x\to \infty} \E (A_x\wedge M) = \E (A\wedge M),\ \ \lim_{M\to \infty} \E (A\wedge M)=\E A\]
and
\[ \lim_{M\to \infty} \E ((B_0-M)1_{\{B_0>M\}}) = 0,\]
we conclude that
\[ \lim_{x\to \infty} \E (A_x) = \E A.\]
\mbox{ } \hfill $\Box$\\

In the next theorem we will discuss the behavior of $X_t$ and
$X'_t$, the second class particle and the dual second class
particle of the process $L_{\la,\mu}$.

\begin{thm}
Let $L_{\la,\mu}$ be the non-stationary Hammersley process defined above. Let $X_t$ be the
$x$-coordinate at time $t$ of a second class particle starting at zero. Then
$$
\lim_{x\rightarrow\infty}P(X_{ax}>x) = \begin{cases} 1,
& a\geq 1/{\mu^2}, \\ {\sqrt{a}-\la}\over{1/\mu-\la}, & 1/{\mu^2}>a\geq \la^2 \\
0,  & \la^2>a\geq 0.\end{cases}
$$
\end{thm}

\noindent {\bf Proof:} To prove this theorem we will calculate the derivative ${\partial\over
\partial h}I(x,h)|_{h=0}$. The proof is inspired by the proof of Theorem 1, \mycite{FerKip}, where
the analogous problem for TASEP is solved.

Fix $t$ and $x$ and consider some small $h>0$. We couple the process $L_{\la,\mu}$ with a process
$\hat L_{\la,\mu}$ constructed in the following way. Consider the same Poisson realizations of
points in the plane. Consider the same realization of sinks. Consider the same realization of
sources on the interval $(h,\infty)\times\{0\}$. On the interval $[0,h]\times \{0\}$ we add an
independent Poisson process of sources of intensity $1/\mu - \la$. The situation is illustrated in
Figure \ref{fig:Proof2_3}.

\begin{figure}[!h]
\begin{center}
\strut
\beginpicture
\footnotesize
\setcoordinatesystem units < 8mm, 6mm>
  \setplotarea x from 0 to 10, y from 0 to 12
\multiput {\bf $\times$} at
4 3
6 5
5 7
/

\put {$0$} at 0.0 -0.5
\put {$x+h$} at 10 -0.5
\put {$t$} at 0.0 8.5
\put {$(x+h,t)$} at 10 8.5
\put {$h$} at 2 -0.5
\put {$x$} at 8 -0.5
\put {$(x,t)$} at 8 8.5
\put {$\bullet$} at 1 0
\put {$X_t^h$} at 10.5 7

\setlinear 
\plot 3 0 3 2 0 2
/
\plot 6 0 6 3 4 3 4 5 0 5
/
\plot 9 0 9 5 6 5 6 6 0 6
/
\plot 10 7 5 7 5 8
/
\plot 0 0 8 0  8 8  0 8 0 0
/
\plot 8 0 10 0
/
\plot 8 8 10 8
/

\setplotsymbol ({$\diamond$}) 
\plotsymbolspacing=8pt
\plot 1 0 1 2 3 2 3 5 4 5 4 6 6 6 6 7 10 7
/

\setplotsymbol ({$\cdot$})
\setdashes <2pt> 
\plot 2 0 2 8 
/

\plot 10 0 10 8
/

\endpicture

\caption{trajectory $X_t^h$ of added point ($\bullet$) hits the right
side and the flux $\bar F_{\la,\mu}(x,t)$ is decreased; for this realization it decreases from 4 to
3.} \label{fig:Proof2_3}
\end{center}
\end{figure}

Note that the process $\hat L_{\la,\mu}$, restricted to the rectangle $[0,h]\times [0,t]$, is a
stationary Hammersley process with intensity of sources $1/\mu$ and intensity of sinks $\mu$. Thus,
the set of crossings of space-time trajectories of $\hat L_{\la,\mu}$ with the vertical segment
$\{h \}\times [0,t]$ is a Poisson process of intensity $\mu$, independent of the sources on
$(h,\infty)\times \{0\}$ and of the Poisson points in $(h,\infty)\times [0,t]$ (this follows from
the stationarity of the Hammersley process when viewed from left to right). Now it follows that the
number of space-time paths of the process $\hat L_{\la,\mu}$ inside the rectangle $[h,x+h]\times
[0,t]$ (which we will denote by $\hat F_{\la,\mu}(x,t)$), has the same distribution as the flux of
the process $L_{\la,\mu}$ on the rectangle $[0,x]\times [0,t]$.

We are interested in the difference
\begin{equation}\label{eq:diffFF}
F_{\la,\mu}(x+h,t)-\hat F_{\la,\mu}(x,t).
\end{equation}

This difference is equal to the difference between the number of space-time paths of the process
$L_{\la,\mu}$ in the rectangle $[0,x+h]\times [0,t]$ and the number of space-time paths of the
process $\hat L_{\la,\mu}$ in the rectangle $[h,x+h]\times [0,t]$.

As we are interested in the expectation of this difference only up to the first order of $h$ when
$h\downarrow 0$, we can distinguish two cases. In the first case there is no extra source in
$[0,h]\times \{0\}$, which means that the space-time paths of $L_{\la,\mu}$ and $\hat L_{\la,\mu}$
are the same. In the second case there is an extra source in $[0,h]\times \{0\}$.

In the first case the difference (\ref{eq:diffFF}) is equal to the number of sources of
$L_{\la,\mu}$ in $[0,h]\times \{0\}$, since those sources are not counted in $\hat
F_{\la,\mu}(x,t)$. Clearly, the expected number of sources in $[0,h]\times \{0\}$ equals $\la h$.

In the second case we have to be more careful. First of all, the probability of having an extra
source equals, in first order, $(1/\mu -\la)h$. The added source can affect the space-time paths of
the original $L_{\la,\mu}$ inside the rectangle $[h,x+h]\times [0,t]$ and can change the flux. Note
that the change in the space-time paths will happen above the trajectory of the added point,
considered as a second-class particle. If this trajectory hits the right side of the rectangle
(i.e. $\{x+h\}\times [0,t]$), it follows that one of the intersections of space-time paths of
$L_{\la,\mu}$ with the right side is eliminated in $\bar L_(\la,\mu)$; see Figure
\ref{fig:Proof2_3}. Thus, in this case the difference (\ref{eq:diffFF}) is $+1$. But if the
trajectory of the added point (which we define by $X^h_t$) does not hit the vertical segment
$\{x+h\}\times [0,t]$, it follows that the sets of intersections are the same for both processes
and the difference is zero again.

Now we can conclude that
\begin{equation}\label{eq:flux}
\E (F_{\la,\mu}(x+h,t)-F_{\la,\mu}(x,t)) = \la h + (1/\mu-\la)h P(X^h_t> x+h) + O(h^2).
\end{equation}

Taking $t = ax$, dividing by $h$ and taking the limit $h\rightarrow 0$ we have

$$
{\partial\over\partial h}I(x,h)|_{h=0} = \la + (1/\mu-\la) P(X_{ax}>x).
$$

From this it follows that ${\partial I(x,h)\over
\partial h}$ is non-increasing in $h\in[-1,1]$ for any fixed $x\geq 1$. Indeed, fix $h_0\in [-1,1]$.
Use (\ref{eq:flux}) to conclude that
\[
{\partial \over
\partial h}I(x,h)|_{h=h_0} = {\partial \over
\partial h}I(x+h_0,h)|_{h=0} = \la + (1/\mu-\la)P(X_{ax}>x+h_0).
\]
The probability $P(X_{ax}>x+h_0)$ obviously does not increase in $h_0$.

We are now in the following situation: we have a family $\{I(x,h):x\geq 1\}$ of functions of $h$,
where $h\in [-1,1]$. For all $x$ the function $I(x,h)$ is differentiable on $[-1,1]$ and concave.
Moreover, we know that the family $I(x,h)$ has a pointwise limit, which is also differentiable on
$[-1,1]$ (Corollary \ref{cor:pwl}). It follows that the family of derivatives also converges to the
derivative of the limit, in other words
$$
\lim\limits_{x\rightarrow\infty}\lim\limits_{h\rightarrow 0} {1\over h }I(x,h)
=\lim\limits_{h\rightarrow 0}\lim\limits_{x\rightarrow\infty}{1\over h}I(x,h).
$$
This is a standard result from analysis.

It follows that
$$
\la + (1/\mu-\la)\lim\limits_{x\rightarrow\infty}P(X_{ax}>x)= \begin{cases} 1/\mu, & a\geq 1/{\mu^2}, \\
\sqrt{a}, & 1/{\mu^2}>a\geq \la^2 \\ \la,  & \la^2>a\geq 0\end{cases}
$$
which finishes the proof of the theorem.\\
\mbox{ } \hfill $\Box$\\

\begin{cor}\label{cor:dual2cl}
Let $L_{\la,\mu}$ be the non-stationary Hammersley process defined above. Let $X'_t$ be the
$x$-coordinate at time $t$ of a dual isolated second class particle at zero. Then

$$
\lim\limits_{x\rightarrow\infty}P(X'_{ax}>x) = \begin{cases} 1,
& a\geq 1/{\mu^2}, \\ {1/\la - \sqrt{1/a}}\over{1/\la-\mu}, & 1/{\mu^2}>a\geq \la^2 \\
0,  & \la^2>a\geq 0.\end{cases}
$$
\end{cor}

\noindent {\bf Proof:} This follows from the symmetry of the trajectories of a second class
particle and a dual second class particle (reversing the role of
$x$ and $t$, interchanging $\la$ and $\mu$ and replacing $a$ by $1/a$).\\
\mbox{ } \hfill $\Box$\\

\section{Second class particles and dual second class particles in the stationary case}

Now we will turn our attention to the stationary case. Let $\la >0$ and consider the Hammersley
process $L_{\la}$ with intensity of sources $\la$, intensity of sinks $1/\la$ and intensity of
$\alpha$-points $1$. As it follows from \mycite{CaGr} (see also Theorem \ref{thm:catgr}), the
trajectories of a second class particle at zero, $(X_t,t)$, and of a dual second class particle,
$(X'_t,t)$, both asymptotically tend to the ray $t = \la^2 x$. We will prove that these
trajectories intersect with probability one and provide some information on the area enclosed by
these trajectories, the origin and the first intersection point.

\begin{thm}\label{thm:touch}
Let $L_{\la}$ be the stationary Hammersley process defined above and let $(X_t,t)$ and $(X'_t,t)$
be the trajectories of a second class and a dual second class particle respectively. Then $(X_t,t)$
and $(X'_t,t)$ intersect with probability one.
\end{thm}

\noindent {\bf Proof:} To prove this theorem we will again use coupling. Consider the standard
coupling of $L_{\la}$ with $L_{\la, 0}$ by removing all sinks. The space-time paths of $L_{\la}$
and $L_{\la,0}$ coincide strictly below the trajectory $(X_t,t)$ (Lemma \ref{lem:leaveout}). Choose a realization such that
$(X'_t,t)$ tends to the ray $t=\la^2x$. If $(X'_t,t)$ and $(X_t,t)$ do not intersect for this
realization, then the trajectory of a dual second class particle of the process $L_{\la,0}$
coincides with $(X'_t,t)$ and, consequently, will also tend to the ray $t = \la^2 x$. Suppose this
event has probability $\eps >0$. Denote $X'_{0}(t)$ as the $x$-coordinate at time $t$ of a dual
second class particle for the process $L_{\la, 0}$. The previous discussion shows that
\[ P(\lim_{x\to \infty } X'_{0}(\la^2x)/x = 1) \geq \eps.\]
Then for any $a>\la ^2$, we would have
\begin{eqnarray}
P(\exists\ M:\forall\ x\geq M: X'_{0}(\la^2x)>\frac{\la^2}{a}x) \geq \eps \nonumber
\end{eqnarray}
This implies that there exists $M>0$ such that
\[P(\forall\ x\geq M: X'_{0}(\la^2x)>\frac{\la^2}{a}x) \geq \eps/2.\]
This in turn implies
\[\forall\ x\geq M: P(X'_{0}(\la^2x)>\frac{\la^2}{a}x) \geq \eps/2\ \  \Leftrightarrow \ \ \forall\ y\geq \frac{a}{\la ^2}M:
P(X'_{0}(ay)>y) \geq \eps/2 \] This contradicts Corollary \ref{cor:dual2cl} for $a$ close enough to
$\la^2$, which proves that $(X_t,t)$ and $(X'_t,t)$ intersect with probability one.\\
\mbox{ } \hfill $\Box$\\

\begin{figure}[!h]
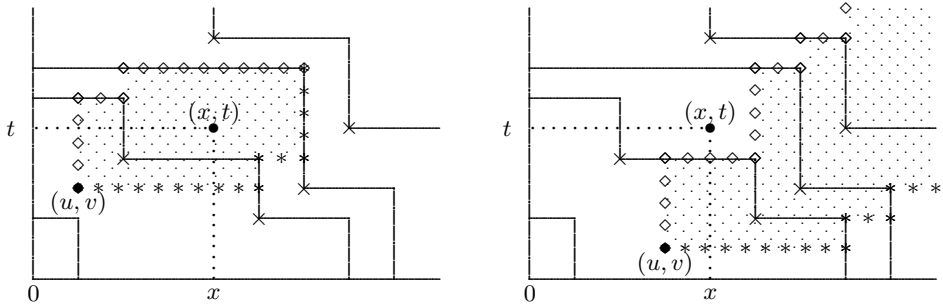

\begin{center}
\strut
\beginpicture
\footnotesize
\setcoordinatesystem units < 6mm, 4mm>
  \setplotarea x from 0 to 25, y from 0 to 12
\multiput {$\times$} at
2 4
5 2
6 3
7 5
4 8
/
   
\multiput {\bf $\bullet$} at
4 5
1 3
/

\put {$0$} at 0.0 -0.5
\put {$(x,t)$} at 4 5.5
\put {$(u,v)$} at 1 2.5
\put {$x$} at 4 -0.5
\put {$t$} at -0.5 5

\setlinear 
\plot 1 0 1 2 0 2
/
\plot 7 0 7 2 5 2 5 4 2 4 2 6 0 6
/
\plot 8 0 8 3 6 3 6 7 0 7
/
\plot 9 5 7 5 7 8 4 8 4 9
/
\plot 0 9  0 0 9 0
/

\setshadegrid span <3pt>

\hshade  3 1 5 4 1 5 
/

\hshade 4 1 6 6 1 6 
/
\hshade 6 2 6 7 2 6
/

\hshade  1 14 18 2 14 18
/

\hshade 2 14 19 3 14 19
/
\hshade 3 14 20 4 14 20
/

\hshade  4 16 20 7 16 20
/

\hshade 7 17 20 8 17 20
/
\hshade 8 18 20 9 18 20
/


\setlinear
\multiput {$\times$} at
13 4
16 2
17 3
18 5
15 8
/
   
\multiput {\bf $\bullet$} at
15 5
14 1
/

\put {$0$} at 11 -0.5
\put {$(x,t)$} at 15 5.5
\put {$(u,v)$} at 14 0.5
\put {$x$} at 15 -0.5
\put {$t$} at 10.5 5

\setlinear 
\plot 12 0 12 2 11 2
/
\plot 18 0 18 2 16 2 16 4 13 4 13 6 11 6
/
\plot 19 0 19 3 17 3 17 7 11 7
/
\plot 20 5 18 5 18 8 15 8 15 9
/
\plot 11 9  11 0 20 0
/

\plotsymbolspacing=7pt
\setplotsymbol ({$*$})

\plot 1 3 5 3 5 4 6 4 6 7
/
\plot 14 1 18 1 18 2 19 2 19 3 20 3
/

\setplotsymbol ({$\diamond$}) 
\plot 1 3 1 6 2 6 2 7 6 7
/
\plot 14 1 14 4 16 4 16 7 17 7 17 8 18 8 18 9
/

\setplotsymbol ({$\cdot$}) 
\setdashes <2pt> 
\plot 15 0 15 5 11 5
/
\plot 4 0 4 5 0 5
/

\endpicture

\caption{In the left picture, the flux at the point $(x,t)$ increases,
for this realization from 2 to 3. In the right picture, the flux at $(x,t)$ remains the same, for
this realization equal to 2.} \label{fig:rainfall}
\end{center}
\end{figure}

Now we will study the expected area enclosed by the trajectories of a second class and a dual
second class particle. The main idea comes from the following observation: suppose we have a
realization of a stationary Hammersley process in the positive quadrant. Now we add an extra
$\alpha$-point at $(u,v)$. How does this change the space-time paths? The situation is depicted in
Figure \ref{fig:rainfall}.

Adding an extra $\alpha$-point at $(u,v)$ moves the first particle immediately to the right of
$(u,v)$ to $(u,v)$; this corresponds to a dual second class particle that starts at $(u,v)$.
However, since there is now a particle at $(u,v)$, the first ``sink'', or crossing on $\{u\}\times
(v,\infty)$, also disappears, which corresponds to a ``normal'' second class particle that starts
at $(u,v)$. Within the region, enclosed by the trajectories of the second class particle and the
dual second class particle, the space-time paths will change in such a way that the flux $F(x,t)$
is increased by 1 for each $(x,t)$ within the trajectories. However, the first time the two
trajectories meet (which happens with probability 1, according to Theorem \ref{thm:touch}), the two
particles annihilate each other (since one corresponds to an extra source, whereas the other
corresponds to an extra sink). This means that the shaded area in Figure \ref{fig:rainfall},
between the two trajectories, is the only area where the space-time paths are changed, and within
this area the flux is increased by 1.

Now we define the following "rainfall" process, which will emulate the adding of extra
$\alpha$-points to a stationary Hammersley process. We define this process in $[0,\infty)^3$, where
the first two variables $x$ and $t$ are the same as we have seen before, but the third parameter
$\tau$ is a kind of external time parameter. Consider three independent Poisson processes in this
space: a Poisson process of intensity 1 on the plane $[0,\infty) \times \{0\}\times [0,\infty)$, a
Poisson process of intensity 1 on the plane $\{0\} \times [0,\infty) \times [0,\infty)$ and a
Poisson process with intensity measure $2\tau dx dt d\tau $ on $[0,\infty)^3$.

Then define $L^{(\tau)}$ as the Hammersley process at ``time'' $\tau$ that uses as its sources the
projection of all Poisson points in $[0,\infty)\times \{0\} \times [0,\tau]$ onto the $x$-axis, as
its sinks the projection of all Poisson points in $\{0\} \times [0,\infty) \times [0,\tau]$ onto
the $t$-axis and as its $\alpha$-points, the projection of all Poisson points in $[0,\infty
)^2\times [0,\tau]$ onto the $(x,t)$-plane. It is easy to check that $L^{(\tau)}$ will be a
stationary Hammersley process with source and sink intensity equal to $\tau$, and the intensity of
$\alpha$-points equal to $\tau^2$. The advantage of this construction lies in the fact that we have
coupled a family of stationary Hammersley processes.

We will use this process to prove the following theorem.

\begin{thm}\label{thm:l2O}
Consider a Hammersley process with intensity of sources and sinks being $\tau$ and intensity of
$\alpha$-points being $\tau^2$. Let $O^{(\tau)}$ be the region in the plane enclosed by the
trajectories $(X_\tau(t),t)$ and $(X'_\tau(t),t)$ of a second class and a dual second class
particle starting in the origin, until their first intersection point. Fix $x,t>0$ and let ${\rm Area}$
denote Lebesgue's measure on the plane. Then
$$
\E {\rm Area}(O^{(\tau)}\cap [0,x]\times [0,t]) = {x+t\over 2\tau} - {1\over 2\tau}(\E (x-X_\tau(t))_+ +
\E (t-X_\tau(x))_+).
$$

\end{thm}
\noindent {\bf Proof:} Consider the rainfall process described above and let $F_\tau(x,t)$ be the
flux at the point $(x,t)$ of the stationary process $L^{(\tau)}$. To prove the theorem we will
calculate the derivative
$$
{\partial\over\partial\tau}\E F_{\tau}(x,t)
$$
in two different ways.

Take some small $d\tau>0$ and note that, since the processes $L^{(\tau)}$ and $L^{(\tau+d\tau)}$
are stationary, we can calculate the expected flux at both times directly:
$$
\E F_{\tau+d\tau}(x,t)-\E F_{\tau}(x,t) = x(\tau+d\tau)+t(\tau+d\tau) - x\tau - t\tau = (x+t)d\tau
$$
and therefore we get
\begin{equation}\label{eq:fluxderdir}
{\partial\over\partial\tau}\E F_{\tau}(x,t) = x+t.
\end{equation}

Now we consider a different way of calculating the difference
\begin{equation}\label{eq:diff}
\E (F_{\tau+d\tau}(x,t)-F_{\tau}(x,t)).
\end{equation}
by using the coupling of $L^{(\tau)}$ and $L^{(\tau +d\tau)}$. Up to the first order of $d\tau$ we
can assume that in the time interval $[\tau,\tau + d\tau]$ there either falls one extra source onto
the segment $[0,x]\times\{0\}$, one extra sink onto the segment $\{0\}\times[0,t]$ or one extra
$\alpha$-point inside the rectangle $(0,x]\times (0,t]$ (if nothing happens the difference is of
course zero).

In the first case the flux will only be changed if the trajectory of a second class particle,
starting at the dropped source on $[0,x]\times \{0\}$, will cross the line $[0,\infty)\times\{t\}$
to the left of the point $(x,t)$. In this case the flux will be increased by one. Similarly, in the
second case the flux will increase by one only if the trajectory of a dual second class particle
starting at the dropped sink on $\{0\}\times [0,t]$ will cross the line $\{x\}\times [0,\infty)$ to
the right of the point $(x,t)$.

The area we are interested in arises from the third case. It follows from the explanation around
Figure \ref{fig:rainfall}, that the flux at $(x,t)$ is increased by one only if $(x,t)\in
O^{(\tau)}(u,v)$, where $(u,v)$ is the location of the extra $\alpha$-point and $O^{(\tau)}(u,v)$
denotes the area described in the theorem, but here the second class particle and the dual second
class particle start at $(u,v)$.

Combining all three cases we have
\begin{eqnarray}\label{eq:firstord1}
\E (F_{\tau+d\tau}(x,t)-F_{\tau}(x,t)) & = & d\tau \int\limits_{0}^{x}P(X_{\tau,t}(u,0)<x)du +
d\tau
\int\limits_{0}^{t}P(X'_{\tau,t}(0,v)>x)dv + \nonumber \\
& & \\
& & +\ 2\tau d\tau \int\limits_0^x\int\limits_0^t P((x,t)\in O^{(\tau)}(u,v))du dv\ +\ o(d\tau),
\nonumber
\end{eqnarray}
where $X_{\tau,t}(u,0)$ and $X'_{\tau,t}(0,v)$ are $x$-coordinates at time $t$ of a second class
and a dual second class particle of $L^{(\tau)}$, starting at points $(u,0)$ and $(0,v)$
respectively.

Using the fact that in the stationary Hammersley process $L^{(\tau)}$, $O^{(\tau)}$ and
$O^{(\tau)}(u,v)$ have the same distribution modulo a translation, we can rewrite the last integral
of (\ref{eq:firstord1}) in the following way:
$$
\int\limits_0^x\int\limits_0^tP((x,t)\in O^{(\tau)}(u,v))du dv =
\int\limits_0^x\int\limits_0^tP((x-u,t-v)\in O^{(\tau)})du dv =
$$
$$
= \E  \int\limits_0^x\int\limits_0^t \mathbf{1}_{O^{(\tau)}}(u,v) dudv = \E {\rm Area}(O^{(\tau)}\cap
[0,x]\times [0,t]).
$$
Dividing both parts of (\ref{eq:firstord1}) by $d\tau$, taking the limit $d\tau\rightarrow 0$ and
using (\ref{eq:fluxderdir}) we arrive at
$$
x+t = \int\limits_{0}^{x}P(X_{\tau,t}(u,0)<x)du + \int\limits_{0}^{t}P(X'_{\tau,t}(0,v)>x)dv +
2\tau \E {\rm Area}(O^{(\tau)}\cap [0,x]\times [0,t]).
$$
Again using the stationarity of $L^{(\tau)}$, we can see that
$$
P(X_{\tau,t}(u,0)<x) = P(X_{\tau}(t)<x-u)
$$
and
$$
P(X'_{\tau,t}(0,v)>x) = P(X'_{\tau}(t-v)>x).
$$
This leads to
\begin{equation}\label{eq:lOans}
\E {\rm Area}(O_{\tau}\cap [0,x]\times [0,t]) = {x+t\over 2\tau} - {1\over
2\tau}\int\limits_{0}^{x}P(X_{\tau}(t)<u)du - {1\over 2\tau}\int\limits_{0}^{t}P(X'_{\tau}(v)>x)dv.
\end{equation}
Consider the first integral on the right side of (\ref{eq:lOans}). We can rewrite it in the
following way:
\begin{eqnarray*}
\int\limits_{0}^{x}P(X_{\tau}(t)<u)du & = & \E (\int\limits_{0}^{x}\mathbf{1}_{\{X_{\tau}(t)<u\}}
du
) \\
& = & \E (x-X_{\tau}(t))_+.
\end{eqnarray*}
Using the equality $P(X'_{\tau}(v)>x) = P(X_{\tau}(x)<v)$, which follows from the symmetry of the
process $L^{(\tau)}$ (note that we have the same intensity of sources and sinks), we conclude a
similar formula for the second integral in (\ref{eq:lOans})
$$
\int\limits_{0}^{t}P(X'_{\tau}(v)>x)dv = \E (t-X_{\tau}(x))_+.
$$

Combining these results we have proved that
\[\E {\rm Area}(O^{(\tau)}\cap [0,x]\times [0,t]) = {x+t\over 2\tau} - {1\over 2\tau}(\E (x-X_{\tau}(t))_+ +
\E (t-X_{\tau}(x))_+).\ \ \ \ \ \ \Box\]

This result is only valid for stationary Hammersley processes with equal source and sink intensity.
However, if we consider $L_\la$, a stationary Hammersley process with source intensity $\la$, sink
intensity $1/\la$ and intensity of $\alpha$-points equal to 1, it is clear that this process has
the same distribution as the image of the process $L_{\la =1}$ under the map
\[ \Phi_\la (x,t) = (x/\la ,\la t).\]
This observation easily leads to the following corollary.
\begin{cor}\label{cor:O_la}
Consider a stationary Hammersley $L_\lambda$. Let $O_\lambda$ be the area in the plane enclosed by
the trajectories $(X_{t},t)$ and $(X'_{t},t)$ of a second class and a dual second class particle
starting in the origin, until their first intersection point. Fix $x,t>0$ and let ${\rm Area}$ denote
Lebesgue's measure on the plane. Then
$$
\E {\rm Area}(O_\la \cap [0,x]\times [0,t]) = {\la x+t/\la \over 2} - {\la \over 2 }(\E (x-X_{t})_+ + \E
(t/\la^2-X_{\la^2 x})_+).
$$
\end{cor}

Another interesting asymptotic corollary to Theorem \ref{thm:l2O} is the following.
\begin{cor}\ Consider a stationary Hammersley process $L_\lambda$.
\begin{enumerate}
\item For any fixed $x_0>0$ we have
$$
\E {\rm Area}(O_{\la}\cap [0,x_0]\times[0,\infty))= {\la x_0+ \la \E X_{\la^2x_0}\over 2}.
$$
\item For any fixed $t_0>0$ we have
$$
\E {\rm Area}(O_{\la}\cap [0,\infty )\times[0,t_0])= {t_0/\la +\la \E X_{t_0}\over 2}.
$$
\item
\[ \lim_{x\to \infty} {\E {\rm Area}(O_{\la}\cap [0,x]\times[0,\la^2 x])\over x} = \la.\]
\end{enumerate}
\end{cor}
\noindent {\bf Proof:} For the first statement, consider $\E (x-X_{t})_+$. Note that
$$
(x-X_{t})_+ \rightarrow 0\ \ \  {\rm a.s.}
$$
as $t\rightarrow \infty$. Since $(x-X_t)_+$ is bounded by $x$, we get
$$
\E (x-X_{t})_+ \rightarrow 0\ \ \ (t\rightarrow \infty).
$$
Now consider $t/\la ^2 - \E (t/\la ^2-X_{\la^2 x})_+ = \E (\min(X_{\la ^2x},t/\la^2))$. Clearly,
$$
\E (\min(X_{\la ^2x},t/\la^2))\rightarrow \E X_{\la^2x}\ \ \ \ (t\rightarrow \infty),
$$
so the first statement is a direct consequence of these observations and Corollary \ref{cor:O_la}.
The proof of the second statement is completely analogous.

The last statement follows from the fact that in $L_\la$, we have
\[ {X_{\la^2 x}\over x} \rightarrow 1\ \ \ \ {\rm a.s.}\]
Now we use dominated convergence to get the desired result from Corollary \ref{cor:O_la}.\\
\mbox{ } \hfill $\Box$\\

This last corollary seems to indicate that the two trajectories of the second class particle and
the dual second class particle run parallel at average vertical distance $\la$. However, the true
behavior of the two trajectories is that they often collide quite quickly, but with low probability
they drift away from each other, in which case $O_\la\cap [0,x]\times [0,\la^2 x]$ grows
quadratically in $x$. It does follow from Theorem \ref{thm:touch} that eventually the two
trajectories will touch. We wish to remark that it is not hard, using the methods of \mycite{CaGr}
and this paper, together with some crude estimates, to show that $\E X_t < +\infty$.

\end{document}